\documentclass{amsart} 
\usepackage{amssymb}
\usepackage{amsmath}
\usepackage{amsfonts}

\begin{document}
\newtheorem{theo}{Theorem}[section]
\newtheorem{prop}[theo]{Proposition}
\newtheorem{lemma}[theo]{Lemma}
\newtheorem{exam}[theo]{Example}
\newtheorem{coro}[theo]{Corollary}
\theoremstyle{definition}
\newtheorem{defi}[theo]{Definition}
\newtheorem{rem}[theo]{Remark}


\newcommand{\Ab}{{\bf A}}
\newcommand{\Bb}{{\bf B}}
\newcommand{\Cb}{{\bf C}}
\newcommand{\Nb}{{\bf N}}
\newcommand{\Qb}{{\bf Q}}
\newcommand{\Rb}{{\bf R}}
\newcommand{\Xb}{{\bf X}}
\newcommand{\Zb}{{\bf Z}}
\newcommand{\Ac}{{\mathcal A}}
\newcommand{\Bc}{{\mathcal B}}
\newcommand{\Cc}{{\mathcal C}}
\newcommand{\Dc}{{\mathcal D}}
\newcommand{\Fc}{{\mathcal F}}
\newcommand{\Hc}{{\mathcal H}}
\newcommand{\Ic}{{\mathcal I}}
\newcommand{\Jc}{{\mathcal J}}
\newcommand{\Lc}{{\mathcal L}}
\newcommand{\Oc}{{\mathcal O}}
\newcommand{\Pc}{{\mathcal P}}
\newcommand{\Sc}{{\mathcal S}}
\newcommand{\Tc}{{\mathcal T}}
\newcommand{\Uc}{{\mathcal U}}
\newcommand{\Vc}{{\mathcal V}}
\newcommand{\Xc}{{\mathcal X}}
\newcommand{\Yc}{{\mathcal Y}}

\newcommand{\ax}{{\rm ax}}
\newcommand{\Acc}{{\rm Acc}}
\newcommand{\Act}{{\rm Act}}
\newcommand{\ded}{{\rm ded}}
\newcommand{\Gm}{{$\Gamma_0$}}
\newcommand{\ID}{{${\rm ID}_1^i(\Oc)$}}
\newcommand{\PAP}{{${\rm PA}(P)$}}
\newcommand{\ACA}{{${\rm ACA}^i$}}
\newcommand{\RefP}{{${\rm Ref}^*({\rm PA}(P))$}}
\newcommand{\RefS}{{${\rm Ref}^*({\rm S}(P))$}}
\newcommand{\Rfn}{{\rm Rfn}}
\newcommand{\tar}{{\rm Tarski}}
\newcommand{\UNFA}{{${\mathcal U}({\rm NFA})$}}

\author{Nik Weaver}

\title [Is set theory indispensable?]
       {Is set theory indispensable?}

\address {Department of Mathematics\\
          Washington University in Saint Louis\\
          Saint Louis, MO 63130}

\email {nweaver@math.wustl.edu}

\date{\em October 26, 2007}

\begin{abstract}
Although Zermelo-Fraenkel set theory (ZFC) is generally accepted as the
appropriate foundation for modern mathematics, proof theorists have known
for decades that virtually all mainstream mathematics can actually be
formalized in much weaker systems which are essentially number-theoretic
in nature. Feferman has observed that this severely undercuts a famous
argument of Quine and Putnam according to which set theoretic
platonism is validated by
the fact that mathematics is ``indispensable'' for some successful scientific
theories (since in fact ZFC is not needed for the mathematics that is
currently used in science).

I extend this critique in three ways: (1) not only is it possible to
formalize core mathematics in these weaker systems, they are in important
ways {\it better} suited to the task than ZFC; (2) an improved analysis of
the proof-theoretic strength of predicative theories shows that most if not
all of the
already rare examples of mainstream theorems whose proofs are currently
thought to require metaphysically substantial set-theoretic principles
actually do not; and (3) set theory itself, as it is actually practiced, is
best understood in formalist, not platonic, terms, so that in a real sense
{\it set theory is not even indispensable for set theory}. I also make the
point that even if ZFC is consistent, there are good reasons to suspect that
some number-theoretic assertions provable in ZFC may be false. This suggests
that set theory should not be considered central to mathematics.
\end{abstract}

\maketitle


Probably most mathematicians are more willing to be platonists about
number theory than about set theory, in the ``truth platonism'' sense that
they firmly believe every sentence of first order number theory has a
definite truth value, but are less certain this is the case for set theory.
Those mathematicians who are unwilling to affirm that the twin primes
conjecture, for example, is objectively true or false are undoubtedly in
a small minority; in contrast, suspicion that questions like the continuum
hypothesis or the existence of measurable cardinals may have no genuine
truth value seems fairly widespread.

Some possible reasons for this difference in attitudes towards number theory
and set theory are (1) a sense that natural numbers are evident and
accessible in a way that arbitrary sets are not; (2) suspicion that sets
are philosophically dubious in a way that numbers are not; (3) the existence
of truly basic set-theoretic questions such as the continuum hypothesis
which are known to be undecidable on the basis of the standard axioms of set
theory, and the absence of comparable cases in number theory; and (4) the
fact that naive set theory is inconsistent. The classical paradoxes of
naive set theory particularly cast doubt on the idea of a well-defined
canonical universe of sets in which all set-theoretic questions have
definite answers.

One philosophically important way in which numbers and sets, as they are
naively understood, differ is that numbers are {\it physically instantiated}
in a way that sets are not. Five apples are an instance of the number 5 and
a pair of shoes is an instance of the number 2, but there is nothing obvious
that we can analogously point to as an {\it instance of}, say, the set
$\{\{\emptyset\}\}$. This is significant because it suggests that we can
be truth platonists about number theory without making any extravagent
metaphysical assumptions about the literal existence of abstract numbers
(``objects platonism''), by interpreting
number-theoretic assertions in terms of these kinds of concrete instances,
whereas it is much harder to see any natural way to meaningfully interpret
set language unless one literally believes in sets as real entities.

Unfortunately, the philosophical difficulties with set-theoretic objects
platonism are extremely severe. First, there is the ontological problem
of saying just what sets are. Second, there is the epistemological
problem of explaining how we, as physical beings, could know anything about
them. And third, there is the fact that naive set-theoretic platonism is
inconsistent and the remedies for this which have been proposed are not
really cogent. (I develop this point in Section 1 below.) Combined with
the observation of the preceding paragraph, this leaves abstract set
theory with no clear philosophical basis.

Perhaps the most influential philosophical defense of set theory is the
Quine-Putnam indispensability argument. According to this argument,
mathematics is {\it indispensable} for various established
scientific theories, and therefore any evidence that confirms these
theories also confirms the received foundation for mathematics, namely
set theory. But as a result of work of many people going back to Hermann
Weyl, we now know that the kind of mathematics that is used in scientific
applications is not inherently set-theoretic, and indeed can be developed
along purely number-theoretic lines. This point has been especially
emphasized by Feferman. Consequently, contrary to Quine and Putnam, the
confirmation of present-day scientific theories provides no special support
for set theory. Indeed, the implications of the indispensability argument
are reversed: if we find that none of the philosophically questionable parts
of mathematics have any clear scientific value, this only tends to strengthen
doubts about their actual validity.

Going beyond the indispensability argument, I see three further questions
which can be raised in defense of the view that axiomatic set theory is the
proper foundation for mathematics. First, even if the bulk of mainstream
mathematics (including all scientifically applicable mathematics) can be
formalized in essentially number-theoretic systems, is the formalization
in Zermelo-Fraenkel set theory (ZFC) not more elegant and aesthetically
appealing? Second, is set theory not needed for {\it some} important
(even if not at present
scientifically applicable) mainstream mathematics? And third,
is set theory not of sufficient intrinsic interest on its own to dissuade
us from abandoning it? I deal with these questions in Sections 4, 5, and
6. Briefly, I claim that formalizing mainstream mathematics in essentially
number-theoretic systems actually carries substantial advantages over
formalization in ZFC; that in fact, essentially {\it all} mainstream
mathematics, not just the scientifically applicable part, is formalizable
in metaphysically unexceptionable systems; and that set theory itself need
not be abandoned because it is already largely practiced in a formalistic
manner
which does not actually require any kind of platonism for its justification.

This raises the possibility that the use of set theory as a foundation for
mathematics may be an historical aberration. We may ultimately find that ZFC
really has no compelling justification and is {\it completely irrelevant} to
ordinary mathematical practice.

\section{Platonist defenses of ZFC}

Before I treat the question of indispensability I want to make a
case that there is no clear philosophical basis for ZFC, on either
platonist or anti-platonist grounds. This section will deal with
platonism and the next with anti-platonism.

According to the platonist view, sets are in some sense real objects
and the Zermelo-Fraenkel axioms are true statements about these objects.
Generally sets are held to be ``abstract'' objects, so that they
are not supposed to exist in space and time. Whatever this even means,
it immediately leads to the {\it epistemological problem} of how we could
know anything about them, a problem given wide currency by Benecerraf
\cite{Ben}.

We must also ask what sort of entities sets are supposed to be (the
{\it ontological problem}). This is related to the epistemological problem;
presumably both questions would have to be answered together. However,
whereas the set theory literature is largely silent on epistemology,
every introductory text has something to say about ontology. Sometimes
this is nothing more than some version of the non-definition ``A set is a
collection of
objects'' (\cite{Mad1}, p.\ 4), but more often the nature of sets is
explained in terms of spurious examples such as flocks or herds. For
instance, Halmos says that a pack of wolves, a bunch of grapes, and a
flock of pigeons are examples of sets (\cite{Hal}, p.\ 1), a statement
satirized by Black with the comment that ``It ought then to make sense,
at least sometimes, to speak of being pursued by a set, or eating
a set, or putting a set to flight'' (\cite{Bla}, p.\ 615).

As Black's comment illustrates, there is little connection between
the mathematical concept of a set and everyday expressions of the kind
cited by Halmos. This was established in a decisive fashion
by Slater, who analyzed in detail the various types of collective
expressions used in ordinary language and showed that none of them
has anything to do with sets in the mathematical sense (\cite{Sla},
Section II). His conclusion is that ``the `sets' of `Set Theory' cannot
be interpreted either in terms of groups of physical things, or in terms
of numbers of things, or by translation into plural expressions. Certainly
there is the set-theoretic symbolism, and the rules for its manipulation,
and maybe it all has some interpretation. But it does not have any of
the traditional interpretations, on the basis of which it was developed''
(\cite{Sla}, p.\ 63). (See also Section 1.1 of \cite{W4} for more on this
point.)

This raises the question of why, if naive set theory really is so
ill-conceived, it is so easy to learn and feels so natural. As I
noted in the introduction, sets are not physically instantiated in
any obvious way. However, they may be at least to some extent
linguistically instantiated,
as natural language exhibits the same hierarchical nested behavior
that is seen in set theory. This suggests
that naive set theory appeals to our lingistic intuitions and that
the classical paradoxes of naive set theory merely reflect limitations
that would be encountered in any attempt to develop an ideal
self-referential language.

Interpreting the paradoxes from a platonist perspective is not so easy.
Indeed, it could be argued that they
should block any solution to the ontological problem, on the grounds that
understanding exactly what sets are is just what would be needed to
render the construction ``the set of all sets'' legitimate.

The standard platonist answer to the paradoxes is that sets should be
understood in terms of the {\it iterative conception} enunciated by
G\"odel: ``a set is something obtainable from the integers (or some
other well-defined objects) by iterated application of the operation
`set of' '' (\cite{God}, p.\ 474). One interpretation of this idea
is that it is meant to single sets out as a special kind of collection.
Some authors clearly display this interpretation by writing, in the context
of discussing the iterative conception, about
collections which do not lie in the cumulative hierarchy and hence
are not sets (e.g.\ \cite{Dra}, p.\ 2 or \cite{Men}, p.\ 40). But this
merely amounts to a change in terminology and evidently leaves the
paradoxes in place as regards collections. Since ``collection'' is
the more primitive concept, this is ultimately unhelpful.

The iterative conception can only help resolve the paradoxes if we
view it as clarifying the intuitive concept of a collection, not as
introducing a new, distinct ``set'' concept. But this clarification
is elusive, as can be seen by looking at some of the versions of the
conception which have appeared in the literature:
\begin{quote}
\noindent sets are `formed', `constructed', or `collected' from their
elements in a succession of stages $\ldots$ (\cite{Par}, p.\ 506)
\end{quote}
\begin{quote}
\noindent According to the iterative conception, sets are created
stage-by-stage, using as their elements only those which have been
created at earlier stages. (\cite{Pot}, p.\ 183)
\end{quote}
\begin{quote}
\noindent In the metaphor of the iterative conception, the steps that
build up sets are ``operations'' of ``gathering together'' sets to
form ``new'' sets. (\cite{Sha}, p.\ 637)
\end{quote}
\begin{quote}
\noindent Thus a set is formed by selecting certain objects $\ldots$
we want to consider a set as an object and thus to allow it to be a
member of another set $\ldots$ When we are forming a set $z$ by choosing
its members, we do not yet have the object $z$, and hence cannot use it
as a member of $z$. (\cite{Sho}, p.\ 322-323)
\end{quote}

It should be apparent from this selection that the nature of the set-forming
operation is extremely unclear. There seems to be no general agreement even
as to whether this is an actual operation which could in any sense be carried
out, or instead some kind of impenetrable metaphor. The problem is apparent
in Boolos's remark that ``a rough statement of the idea $\ldots$ contains
such expressions as `stage', `is formed at', `earlier than', `keep on going',
which must be exorcised from any formal theory of sets. From the rough
description it sounds as if sets were continually being created, which is
not the case'' (\cite{Boo}, p.\ 491). Yet Boolos does not follow his rough
statement with a more informative informal description that avoids the
objectionable phrases, and it seems doubtful that he could. Without these
expressions there {\it is} no informal description.

This difficulty is connected to the ontological problem, about which none of
the authors cited above has anything meaningful to say: if we have no idea
what sets are supposed to be, obviously there is little we can say about
how they are supposed to be formed. Yet the idea that sets are in some
sense ``formed'' from elements which enjoy some kind of ``prior'' existence
is crucial to the iterative conception's ability to evade the classical
paradoxes. The point is supposed to be that the set of all sets, or the
set of all ordinals, or Russell's set, are illegitimate on the iterative
account {\it precisely
because they cannot be ``formed''}. So it seems fair to say that the
iterative conception successfully deals with the paradoxes only to the
extent that it presents a clear picture of the operation of set formation,
which is to say, not at all.

This problem is especially acute because it is also part of the
iterative conception that the power set operation is the basic step to
be used in the ``construction'' of sets. Jan\'e explains why this causes
trouble (a similar point is made by Lear \cite{Lea}):
\begin{quote}
We must also have recourse to some form of the power set operation
before setting up the iterative conception. This is an important point
that is often obscured and whose neglect might lure us into believing
that the power set axiom of ZF simply follows from the idea of iteration.
The reason given for the validity of this axiom is that if a set lies on
a layer, so do all its subsets, and therefore the set of all of them lies
on the next layer. One question about this way of presenting the matter
is what is meant by ``all subsets'' of a set $a$. Perhaps from the
standpoint that the iterative conception only describes how the world
of sets is actually structured there is really no question to be asked
(for if we can resort to the universe of sets, there is no difficulty
in saying what are all subsets of $a$; they are just those sets all of
whose members are members of $a$). But if we want to account for the
set-theoretic universe as built by iterated application of the power
set operation, such an explanation is of no use whatever. Since we
cannot turn to the result of the iteration to tell what to do at each
step, the notion of {\it all subsets} of a given set cannot be taken
for granted, but must be clarified at the outset. (\cite{Jan}, p.\ 374)
\end{quote}

This actually understates the problem, because however the notion of
{\it all subsets} is clarified, if the force of the iterative conception
against the paradoxes is to be maintained this must be done in a manner
that allows us to retain some sense of ``construction'' or ``formation''
of power sets. But the essentially circular nature of power sets renders
this prospect quite hopeless. I refer here to the predicativist criticism
of the power set of the natural numbers as an inherently circular object.
The sort of circularity that is involved can be seen in an example
I introduced in \cite{W1}. Let $(A_n)$ be a standard enumeration of the
sentences of second order arithmetic. These are formal expressions which
include first order variables $x, y, z, \ldots$ ranging over the natural
numbers and second order variables $X, Y, Z, \ldots$ ranging over sets of
natural numbers, and which allow quantification of both types of variables.
Then consider the set $S = \{n: A_n$ is true$\}$. This is a well-defined
set of numbers, provided each sentence $A_n$ has a definite truth value.
However, it is the existence of the power set of the natural numbers as a
well-defined totality which guarantees that these sentences do have definite
truth values; without assuming this we cannot give a definite meaning to
second order quantification. Thus, if the power set of the natural
numbers is not already available then we cannot expect to be able to
determine the set $S$, or to put it differently, {\it $S$ is a set of
natural numbers which only becomes available after we have all sets of
natural numbers}. Because of this kind of circularity, there is no sense
in which we can imagine the power set of the natural numbers as being
built up in a piecemeal fashion, if it was not already available
to begin with.

The ``medieval metaphysics'' (\cite{Fef}, p.\ 248) of platonism must be
attacked directly on epistemological and ontological grounds. But even
leaving these questions aside, the iterative conception is internally
incoherent. If language about ``formation'' and ``construction'' is to be
taken as having any content whatever --- which it must if the paradoxes
are to be defused --- then the idea that the formation of power
sets is to count as a legitimate set-building operation is meretricious.
There is no meaningful sense in which power sets can be thought of as
being formed or constructed.

\section{Anti-platonist defenses of ZFC}

The goal of this section is to argue that there is no good anti-platonist
justification of ZFC. But first something needs to be said about what
would count as a ``justification'' in this context. Since anti-platonists
do not believe in a well-defined universe of sets, they obviously cannot
be expected to affirm that the Zermelo-Fraenkel axioms are really true in
any strong sense.

It is tempting to suppose that the only issue for anti-platonists is
whether ZFC is consistent. However, this ignores the fact that {\it some}
theorems of ZFC are genuinely meaningful even from a very strong
anti-platonist perspective. For example, certain sentences in the
language of set theory are directly interpretable as statements of
first order arithmetic. Anyone who recognizes such statements as
meaningful ought to be interested in knowing whether all theorems
of ZFC of this form are actually true.

At a minimum, I think that any attempt to justify ZFC on anti-platonist
grounds ought to say something about whether ZFC is {\it $\Sigma_1$-valid},
i.e., whether every $\Sigma_1$ sentence of first order arithmetic that is
provable in ZFC is actually true.
These are, in effect, assertions of the form ``Turing machine $x$ halts
on input $y$'' for particular values of $x$ and $y$. I take it that few
anti-platonists would deny the meaningfulness of such assertions. Therefore
the question of the $\Sigma_1$-validity of ZFC is legitimate. It is also
pressing: if ZFC gives us bad information about which Turing machines halt,
we surely would not want to use it as a basis for mathematical theories
used in scientific applications, nor should we be enthusiastic about using
it as a basis for mathematical theories generally.

Mere $\Sigma_1$-validity would seem to be a very minimal requirement,
but there are a number of purported anti-platonist justifications of ZFC
which in fact do not address this issue in any way. For example, a variety
of reasons have been given by various authors
as to why anti-platonists ought to feel confident that
ZFC is consistent. These arguments may or may not be persuasive, but even
if we could be certain that ZFC is consistent, that in itself would not
entail its $\Sigma_1$-validity. For instance, if ZFC
is consistent then the theory ${\rm ZFC} + \neg\,{\rm Con}({\rm ZFC})$ is
also consistent, but it proves the false $\Sigma_1$ sentence
$\neg\,{\rm Con}({\rm ZFC})$.

(According to (\cite{Smo}, pp.\ 1-2), G\"odel criticized the formalist goal
of proving the consistency of infinitary mathematics in precisely this way:
``\ldots it would be, e.g., entirely possible that one could prove with
the transfinite methods of classical mathematics a sentence of the form
$\exists x F(x)$ where $F$ is a finite property of natural numbers (e.g.\
the negation of the Goldbach conjecture has this form) and on the other
hand recognise through conceptual considerations that all numbers have
the property {\it not}-$F$; and what I want to indicate is that this
remains possible even if one had verified the consistency of the formal
system of classical mathematics.'')

``Naturalistic'' philosophies which rationalize ZFC in terms of the fact
that set theorists like its properties also fail the $\Sigma_1$-validity
test. For example, Maddy approves of arguing from ``this theory has properties
we like'' to ``this theory is true'' (\cite{Mad2}, p.\ 163); but we may like
a theory for reasons which have no bearing on its $\Sigma_1$-validity,
and the sorts of reasons Maddy admires generally do fall in this category.
In particular, she rejects demands that mathematical theories be true
in any substantive sense as ``philosophical niceties'' and states that
``what matters are the intra-mathematical goals and the effectiveness of
various means of achieving them'' (\cite{Mad3}, p.\ 417). Probably when
she says this she is thinking of questions like the continuum hypothesis
which could plausibly be supposed to have no definite content, and not
about the fact that axioms which settle such questions might also have
elementary number-theoretic consequences.

On the other hand, there are more ambitious philosophical programs, such as
those in \cite{Chi} and \cite{Fie}, which attempt to justify ZFC (or a large
fragment of ZFC) in a very strong way on anti-platonist grounds. If some
program of this sort were to succeed, we could then be confident that ZFC
is $\Sigma_1$-valid. However, these attempts are, I believe, generally
regarded as
actually involving substantial platonist assumptions, and this is only to
be expected given the proof-theoretic strength of ZFC. As I suggested in
the introduction, a straightforward case can be made for Peano arithmetic
on anti-platonist grounds, but when we turn to set theory the prima facie
case goes against anti-platonist acceptability. Indeed, Peano arithmetic
(PA) lies more or
less at the limit of what can be regarded as obviously legitimate by an
anti-platonist. There are accepted methods of working up to stronger
systems, for example by adding an assertion of the consistency of PA and
then iterating (cf.\ Section 5), but it seems extremely unlikely that one
could get all the way up to ZFC using such techniques. The point is that
there is a {\it huge} gap in consistency strength between PA and ZFC, so
that the prospect of incrementally working up to ZFC appears quite
hopeless. The other possibility, that there is some anti-platonist principle
that would allow us to bridge the gap in a single step, also seems highly
unlikely. So the prospect of anti-platonistically justifying ZFC in a
manner sufficient to establish its $\Sigma_1$-validity is not realistic.

I have been taking it as prima facie plausible that Peano arithmetic is
anti-platonistically legitimate; I defend this claim in more detail in
Section 1.3 of \cite{W5}. Various objections could be raised. For example,
some set-theoretic platonists would say that even number theory requires
objects platonism. Their response to the argument that numbers can be
understood in terms of concrete instances might be that this only works
for small numbers since the universe could be finite. To the contrary,
I find it perfectly reasonable to suppose that one can understand basic
arithmetic without having to believe in an abstract world of numbers.
The question whether our universe is finite or infinite does not seem
relevant since the point is not that we have to actually physically
observe and manipulate $n$ objects in order to have access to the number
$n$. All we need is to have an idea of what it would mean to observe and
manipulate $n$ objects. I think this is a perfectly sensible way to
interpret number theory and I think the idea that there is something
deeply problematic with it is disingenuous. Note that set theory
apparently cannot be interpreted in any similar way because there are
no natural concrete proxies for sets.

An objection can also be made in the other direction, from the point
of view of intuitionism. Whereas the set theorist accepts PA himself,
but argues that one must be an objects
platonist to do this, the intuitionist is genuinely unable to accept PA.
This position has integrity but I think it runs against common sense.
Even if we reject the idea of a platonic world of numbers and interpret
number-theoretic assertions in concrete terms, questions like the twin
primes conjecture still seem to me completely definite. I argue this
point further in Section 1.3 of \cite{W5}.

It is important to recognize that the legitimacy of ZFC as a foundational
system involves more than its mere alleged consistency. If it is to be
taken as the rightful foundation of mathematics we should at the very least
have good reasons for believing it is $\Sigma_1$-valid. Yet ``soft''
justifications of ZFC which merely aim to show that it is consistent
or pleasing in some way do not speak to this question. There remains the
possibility of a ``hard'' justification which does accomplish something
substantive on this score, but this does not seem a realistic hope.

I will go further and say that it is more likely than not that some false
statements of first order arithmetic are theorems of ZFC; see Section 7
below.

\section{Applicable mathematics without set theory}

In the preceding two sections I have indicated that Zermelo-Fraenkel
set theory does not have a clear philosophical basis in either platonist
or anti-platonist terms. It therefore becomes reasonable to ask what the
consequences would be of rejecting ZFC as a foundation for mathematics.

Some philosophers may naturally be reluctant to pursue this question
because it could entail having to tell mathematicians that they are
practicing their subject incorrectly. The situation is not quite as bad
as that, since, after all, most mathematicians have little interest in
foundations and may have no particular commitment to ZFC. (Maddy \cite{Mad2}
paints a very different picture, but her ``mathematicians'' really seem
to be set theorists.) Still, for example, one
commentator, who clearly recognized the difficulty in justifying
the power set axiom, was driven to simply postulate the existence of power
sets based on ``external requirements'', presumably meaning the fact that
mathematicians use them (\cite{Jan}, p.\ 388).

However, it is now a settled fact that power sets of infinite sets are not
actually needed for the vast bulk of mainstream mathematics. The philosophical
stance which admits the natural numbers but not its power set is called
{\it predicativism}; it was originally put forward by Bertrand Russell
and Henri Poincar\'e, and there is a long line of research stretching back
to Hermann Weyl which establishes in detail its ability to encompass
ordinary mathematics. (See the introduction to \cite{W2} for references,
and see \cite{W1} and Section 3 of \cite{W4} for more on the philosophical
basis of this view.) The basic idea is that we accept the natural numbers
and individual real numbers (or equivalently, individual sets of natural
numbers, which can still be pictured in terms of physical instantiation) but
do not assume the existence of a well-defined set of all real numbers (which
cannot be meaningfully understood in terms of physical instantiation). In
effect we treat the real line as a proper class. We can then
accomodate all of the structures that appear in normal mathematics by
various encoding tricks; for example, using an injection from
$\Nb \times \Nb$ into $\Nb$
we can encode a sequence of reals as a single real number,
and so on. Formalizing this approach yields the system ${\rm ACA}_0$, and
the fact that the bulk of normal mathematics can actually be carried out
in this system is established at length in \cite{Sim}. ${\rm ACA}_0$ is
conservative over PA, i.e., any theorem of first order arithmetic provable
in ${\rm ACA}_0$ is already provable in PA.

There are a number of possible variations on the above. For instance,
in \cite{W5} I present a system which allows a limited additional
ability to reason about arbitrary properties of real numbers. This has
to be done carefully in order to maintain the philosophical integrity
of the system, but it has the advantage of reducing the coding
machinery needed in \cite{Sim}. Another approach is to explicitly
construct a countable structure which can play the role of a miniature
set-theoretic universe, and then formalize core mathematics within this
structure in the standard way. This was done in \cite{W2}.

Feferman has made the point that the fact that normal mathematics can
be formalized in systems like ${\rm ACA}_0$ neutralizes the Quine-Putnam
indispensability argument mentioned in the introduction (\cite{Fef},
Chapter 14). The latter is possibly the most influential philosophical
defense of set-theoretic platonism, but in light of the above facts it
completely loses its force. If all scientifically applicable mathematics
can be straightforwardly formalized in systems which are conservative over
PA, then the fact that some scientific theories which use this mathematics
are well confirmed lends no particular support to ZFC.

Set theory apologists have responded rather grudgingly to Feferman's point.
Quine himself, writing in 1991, acknowledged its validity in a magnanimous
but curiously subjunctive way:
\begin{quote}
Sanguine spirits there have been, and Solomon Feferman and Hao Wang
are two of them today, who hope to show that enough mathematics can
be derived for purposes of natural science without going beyond
predicative set theory. This would be a momentous result. It would
make a clean sweep of the indenumerable infinites and unspecifiable
sets $\ldots$ (\cite{Qui}, p.\ 229)
\end{quote}
The ``hope to show'' and ``would'' language is a little strange, given
that Wang outlined how to do this in 1954 \cite{Wan1} and gave details in
1963 (\cite{Wan2}, Chapter XXIV), and several others had done similar
work between 1955 and 1975.

Tellingly, Maddy finds it ``striking'' (\cite{Mad1},
p.\ 4) and ``remarkable'' (\cite{Mad2}, p.\ 23) that so much mathematics
can be formalized in ZFC, and concludes ``that set theory plays this
role is central to modern mathematics'' (\cite{Mad2}, p.\ 35), yet as
far as predicative systems are concerned she makes only the vague
admission that ``for the purposes of providing tools for current science
$\ldots$ weaker systems would probably do'' (\cite{Mad3}, p.\ 413) and
sees no special significance in this fact. Steel also considers the
formalizability of mathematics in ZFC ``remarkable'' (\cite{Ste},
p.\ 423) but thinks the fact that mainstream mathematics can be done in
(unspecified) weaker systems ``will never be more than a description of the
current state of affairs'' (\cite{Ste}, p.\ 424). In both cases formalizability
in the favored system, ZFC, is accorded great significance but the stronger
phenomenon of formalizability in essentially number-theoretic systems is
dismissed as irrelevant.

Hellman actually claims that at present ``{\it some} physically applicable
mathematics appears to transcend the bounds of predicativism, especially
the use of nonseparable Banach or Hilbert spaces, e.g.\ in quantum field
theory (see e.g.\ Emch (1972))'' (\cite{Hel}, p.\ 218; italics in
original). It is unfortunate
that this reference is not more specific, as I cannot find any use of
nonseparable Hilbert spaces in Emch's book. They certainly play no
significant role in modern quantum field theory, or in any area of modern
physics for that matter.

Nonseparable Banach spaces, on the other hand, are routinely used in the
mathematics of quantum field theory: every infinite dimensional von Neumann
algebra is nonseparable. But already in ordinary quantum mechanics, the space
$B(H)$ with $H$ separable is itself non-separable,
and even in engineering applications nonseparable $L^\infty$ and $H^\infty$
spaces are commonly used, so one hardly has to go all the way to quantum
field theory to see nonseparable Banach spaces in applications. However,
every such example of which I am aware (including local nets of von
Neumann algebras in quantum field theory, presumably the point of the
reference to Emch) can be straightforwardly handled in, e.g., the
systems of \cite{W2} or \cite {W5} or, with minor indirection, in
${\rm ACA}_0$. (The Banach spaces cited above are all
weak* separable, local nets of von Neumann algebras are countable
unions of such spaces, and so on.) So Hellman's assertion is not
well-taken.

Hellman, Maddy, and Steel are all impressed by the possibility that
set-theoretically substantial mathematics might one day be needed in
scientific applications. But of course the mere possibility
of future applications provides no support whatever for the indispensability
argument. Indeed, one could say of virtually any formal system that future
applications are {\it possible}. We would only find this possibility
noteworthy if we had separate reasons for being interested in the particular
system in question. This is really the opposite of an indispensability
argument, because ZFC is not gaining credibility from its scientific
applications --- at present it has none --- but rather is seen as a
good candidate for future applications because evidently it is
already felt to be credible on some {\it other} grounds.

I must add, however, that given our current understanding of basic physics,
the prospect of set-theoretically substantial mathematics ever becoming
essential to
meaningful scientific applications appears {\it extremely} unlikely.
This should be obvious to anyone with a basic knowledge of mathematical
physics and an understanding of the scope of predicative mathematics. An
essential incorporation of impredicative mathematics in basic physics would
involve a revolutionary shift in our understanding of physical reality of a
magnitude which would dwarf the passage from classical to quantum mechanics
(after all, both of these theories are completely predicative). I would rate
the likelihood of ZFC turning out to be inconsistent as much higher than the
likelihood of it turning out to be essential to basic physics.

The assumption that set-theoretically substantial mathematics is of
any use in current science is simply false.
One can hold out hope that some radically different future physical
theory would require such mathematics, but there is no rational basis
for this hope. It certainly finds no encouragement in the character
of current physics. The argument that ZFC should be retained as the
standard foundation for mathematics because it might conceivably be
indispensible to some future scientific application only makes sense
if we have an independent reason for favoring ZFC over other
foundational systems, and is not itself a reason for favoring ZFC.
Consequently, given the scope of predicative mathematics, the Quine-Putnam
indispensability argument no longer has any force whatever.

\section{ZFC is foundationally unsuitable}

The idea that Zermelo-Fraenkel set theory is to be justified by its
scientific applications is not the only potentially persuasive argument
in favor of taking it to be the appropriate foundational system for modern
mathematics. A case could also be made that predicative systems, while in
principle adequate, are too awkward in practice and should yield to ZFC
on aesthetic grounds. Or one could argue that in order to be acceptable a
foundational stance must encompass all mainstream mathematics, not just the
scientifically applicable parts. Or one could argue directly in defense
of ZFC as a mathematical system of intrinsic interest.

I will address these points in this and the following two sections. First,
is ZFC really the most attractive foundational system available? It is
certainly aesthetically preferable to some predicative systems, particularly
some of the older ones, but that does not settle the issue. A variety of
predicative systems have been proposed as settings in which to formalize
core mathematics; the question is whether any of them competes with ZFC
in terms of elegance and ease of use.

For instance, several predicative systems based on type-theoretic
formalisms have been put forward by various authors. Personally, I
tend to dislike this sort of approach. Partly this is because extra
work seems to be involved in keeping track of the different types, and
partly it is because I find some of these systems unintuitive, but
probably my main disagreement with type-theoretic approaches generally
is that they seem stylistically too far removed from mainstream
mathematical practice. Probably this is simply a matter of taste.

A different complaint can be raised against the formalization of core
mathematics in ${\rm ACA}_0$ as described in \cite{Sim}, namely that
it involves fairly heavy coding machinery. However, that criticism is
unfair because the goal in \cite{Sim} is to determine the {\it weakest
possible} systems in which various theorems can be proven. Getting by
with absolutely minimal assumptions may require some extra coding, but
that is not the point. As I show in \cite{W5}, passing to a third order
language substantially reduces the need for coding.

We must remember that every formalization of mathematics involves
some sort of coding. In ZFC natural numbers are coded as von Neumann
ordinals, integers are coded as equivalence classes of ordered pairs
of natural numbers (and ordered pairs are coded set-theoretically),
rationals are coded as equivalence classes of ordered pairs of integers,
and reals are coded as (say) Dedekind cuts. The predicative system CM
of \cite{W5} requires roughly the same degree of coding. Here natural
numbers are taken as primitive, so we do not have to code them as von
Neumann ordinals. We use an injection from $\Nb \times \Nb$ into $\Nb$
to code ordered pairs of natural numbers as natural
numbers; I do not think this is terribly awkward compared to coding
ordered pairs set-theoretically. The constructions of $\Zb$, $\Qb$, and
$\Rb$ in CM are then similar to their constructions in ZFC.

With one additional use of an injection from $\Nb \times \Nb$ into $\Nb$
the preceding constructions yield sequences of integers,
sequences of rationals, and sequences of reals. Now virtually every standard
mathematical space can be more or less straightforwardly realized inside
$\Rb^\omega$, so the amount of coding needed in CM to construct these
standard spaces is roughly comparable to that needed in ZFC.

One of the less attractive coding aspects of ${\rm ACA}_0$ involves
continuous maps between complete separable metric spaces. In this
setting a complete separable metric space $X$ is coded by a countable
dense subset of $X$, and functions
between such spaces are slightly complicated to handle because we
cannot assume that the given dense subset of the domain is mapped into
the given dense subset of the range. In CM this difficulty disappears
because the third order language allows us to represent separable spaces
(and even some nonseparable spaces) directly.

In short, I think the formalization of core mathematics in CM is quite
comparable with its formalization in ZFC in terms of simplicity,
elegance, and ease of use. Of course real mathematics is not actually
formally executed in ZFC; it is presented informally in a manner which,
ideally, would render formalization in ZFC tedious but not difficult.
Instead taking CM as the foundational standard would in most subjects
alter this informal presentation of mathematics in everyday practice
{\it not at all}. In some more set-theoretically oriented fields like
functional analysis there would be noticeable differences, but still,
mainstream practice would not have to change in any really substantive
way. The idea that predicative mathematics has to be horribly
complicated is just not true.

The formalizations of core mathematics in ZFC and CM are roughly comparable
in terms of elegance. However, this does not mean that the two systems should
be thought equally suitable for this task. ZFC has one major shortcoming that
predicative systems do not share, namely, that {\it the Zermelo-Fraenkel
universe is grossly discordant with the realm of ordinary mathematics.}
Following G\"odel's iterative conception of sets (see Section 1), we can
prove in ZFC the existence
of sets $S_\alpha$, with $\alpha$ ranging over all ordinals, such that
$S_0 = \Nb$, $S_{\alpha+1} = \Pc(S_\alpha)$ (the power set of $S_\alpha$),
and $S_\alpha = \bigcup_{\beta < \alpha} S_\beta$ when $\alpha$ is a
limit ordinal. This means that $S_\alpha$ is defined for $\alpha = \omega$,
$\omega^2$, $\omega^\omega$, $\epsilon_0$, $\aleph_1$, $\aleph_\omega$,
$\aleph_{\aleph_\omega}$, and so on. Yet {\it virtually none of this
sequence beyond $S_1 = \Pc(\Nb)$ is needed in mainstream mathematics}.
Almost no objects, arguably no objects at all, in mainstream mathematics
have cardinality greater than that of the continuum, and consequently
virtually every ordinary mainstream object can be more or less
straightforwardly encoded as either a set of natural numbers or a set
of reals. This dichotomy between sets of numbers and sets of reals
is just the dichotomy between {\it discrete} and {\it continuous}
mathematics. We have no analogous word for mathematics at the level of
$\Pc(\Pc(\Nb))$ or at any higher level because there is no mainstream
mathematics there. Maybe the right word is {\it pathological}.

In contrast, the predicative universe exhibits a strikingly exact fit with
the universe of ordinary mathematics. This is particularly well illustrated
by the dual Banach space construction in functional analysis. Classically,
every Banach space $V$ has a dual Banach space $V'$, but predicatively
this construction is only possible when $V$ is separable. Remarkably,
it seems to be a general phenomenon that for any ``standard'' Banach
space $V$, its dual is also ``standard'' if and only if $V$ is separable.
In other words, start with any well-known Banach space $V$ that commonly
appears in the functional analysis literature, and iteratively take duals
to create a sequence $V, V', V'', \ldots$. It will generally be the case
that if $V$ is in common use then $V'$, $V''$, etc., will also be in
common use --- up to the first nonseparable space in the sequence.
All spaces after that point will be highly obscure.

For example,
take $V = L^1(\Rb)$ (separable). Then $V' = L^\infty(\Rb)$ (nonseparable),
and $V''$ is an obscure space that has no standard notation. Or take $V =
C[0,1]$ (separable). Then $V' = M[0,1]$ (nonseparable), and $V''$ is
an obscure space that has no standard notation. Take $V = K(\Hc)$ (the
compact operators on a separable Hilbert space $\Hc$; separable). Then
$V' = TC(\Hc)$ (the trace class operators on $\Hc$; also separable),
$V'' = B(\Hc)$ (the bounded operators on $\Hc$; nonseparable), and
$V'''$ is an obscure space that has no standard notation. If $V = c_0$
(separable) then $V' = l^1$ (separable), $V'' = l^\infty$ (nonseparable),
and $V'''$ arguably has a standard notation --- it is the space of Borel
measures on the Stone-\v{C}ech compactification of the natural numbers ---
but it is certainly an obscure space that appears in the literature with
extreme rarity. Examples of this type could be multiplied endlessly.

The explanation of this phenomenon is simple: generally speaking, duals
of nonseparable
spaces are highly pathological objects about which little of value
can be said. This is characteristic of impredicative mathematics
generally. The extra generality of ZFC is spurious, involving structures
which are highly pathological from the point of view of mainstream
mathematics. Generally speaking, ``nice'' spaces are predicative. This
may sound like
a purely subjective judgement, but it can also be seen more objectively
in the fact that basic properties of impredicative spaces tend to be
undetermined in ZFC. For instance, the most familiar example of a
pathological, impredicative space is $\beta\Nb$, the Stone-\v{C}ech
compactification of the natural numbers. Does $\beta\Nb - \Nb$ have any
nontrivial self-homeomorphisms? The answer is independent of ZFC, assuming
ZFC is consistent \cite{vM}. The existence of basic questions which
cannot be answered in ZFC is typical of impredicative spaces.

Of course, the existence of undecidable statements is a feature of any
sufficiently strong consistent formal system. The issue is not whether
such statements exist, but how common and how basic they are.
The impredicative portion
of the universe of ZFC is rife with undecidability, and many known
undecidable statements are not contrived, but appear quite
fundamental --- the prototype of such a statement being the continuum
hypothesis. Despite strenuous efforts by logicians to identify similarly
fundamental number-theoretic statements which are undecidable in
predicative systems, the best available examples are still rather
complicated and unnatural. Extravagent claims along these lines are
sometimes made, but this seems to be wishful thinking similar to the
expectation that impredicative mathematics is likely to become
scientifically relevant (see the end of Section 3). Even if some really
good examples were found, this still would not compare with the ubiquity
of undecidability in impredicative mathematics.

The pathological quality of ZFC has consequences. One effect is
that working mathematicians in certain fields have to expend effort
learning to avoid set-theoretically pathological lines of investigation.
This is illustrated by Hellman's suggestion mentioned in Section 3
that nonseparable Hilbert spaces could be important in mathematical
physics. They are not, but this is something mathematical physicists
have to learn. Similarly, functional analysts have to learn that
duals of nonseparable spaces are not fruitful. Thus, in some fields the
greater generality of ZFC merely opens non-productive avenues that working
mathematicians must actively avoid.

I do not want to overemphasize this point, as it is not difficult to
recognize set-theoretic pathology or to learn that it is to be avoided.
However, the standard use of ZFC as a foundation for mathematics does
entail some wasted effort in this way. A more significant effect is
the loss of top-flight mathematicians to impredicative set theory. This is
not to say that the study of ZFC should be abandoned wholesale (see Section
6), but on the other hand it is obvious that the intellectual power that
has gone into the development of set theory is far out of proportion to its
importance to mathematics as a whole. The erroneous assessment of ZFC as
being central to mathematics has attracted many first-rate mathematicians
to study it, drawing them away from other genuinely central subjects.

ZFC is not an appropriate foundation for mathematics. While the formalization
of core mathematics in ZFC is reasonably straightforward, it is no more
elegant than formalization in predicative systems like CM. But ZFC, unlike
CM, is very poorly fitted to mainstream mathematics in that it embeds the
well-behaved realm of ordinary mathematics in a vast arena of set-theoretic
pathology. This can be a distraction for ordinary mathematicians because it
opens up fruitless lines of investigation. More significantly, the widespread
ascription of fundamental status to ZFC has the effect of channeling
intellectual resources away from truly central subjects. Taking ZFC as
the foundational standard is, in important ways, pernicious.

\section{The limits of predicativity}

As I discussed in Section 3, current science makes no essential use
of impredicative mathematics. One reaction that some people have to
this fact is that it is just an accident and future
scientific theories surely will (or at least, there is a reasonable
expectation that they will) use impredicative mathematics. I think
this could seem plausible if one thinks of mathematics as a
whole as being highly interconnected in an intellectually and
aesthetically compelling way, and one also sees predicative
mathematics as a small, artifically restricted part of mathematics
as a whole. Then one might see no significance in the fact that all
currently applied mathematics happens to be predicative and no
reason to expect this to remain true as science develops further.

The clearest error in this way of thinking is the idea that predicative
mathematics is a {\it small} part of mathematics as a whole. On the
contrary, most mainstream subjects --- differential geometry, algebraic
topology, complex analysis, PDEs, etc. --- as they are currently
practiced, lie virtually entirely within the bounds of predicativity. Some
other mainstream areas like abstract algebra or functional analysis in
principle include impredicative material, but the role impredicativity
plays in current research in these fields is still quite minimal. It is
only in set theory itself that significant impredicativity is routinely seen.

This calls into question the ``organic unity'' picture of mathematics as
regards set theory. The latter is {\it not} interconnected with the other
fields listed above in anything like the same way that those fields are
interconnected with each other. This should be obvious to any working
mathematician; it can also be seen, for example, in the fact that virtually
all mainstream mathematics can be straightforwardly construed as taking
place in (at worst) $\Pc(\Pc(\Nb))$.

Thus, Feferman's emphasis on the fact that scientifically applicable
mathematics is predicative (see Section 3) substantially understates
the scope of predicativity, as it can leave the impression that there
may be a large body of scientifically unapplied mainstream mathematics
which is not predicative. Certainly, framing the debate in terms of
scientifically applicable mathematics is the best way to make a case
specifically against the Quine-Putnam indispensability argument because
it is only this type of mathematics that is relevant to that argument.
However, it may not be persuasive as a general defense of predicativism
because of course one wants to preserve all mainstream mathematics, not
just those bits that are currently being used in science.

Now it is more difficult to decisively show that all or virtually
all mainstream mathematics is predicative because ``mainstream'' is not
as sharp a concept as ``scientifically applicable''. But there are still
ways of making this case. For example, one can observe that the material
that appears on qualifying exams in typical mathematics Ph.D.\ programs is
entirely or almost entirely predicative. By making observations like this
one can build up the idea that impredicative mathematics is, by mainstream
standards, highly pathological, and this should lead to the expectation
that it is unlikely to be significantly used in successful
scientific theories.

Although it is imprecise, the question of just how much mainstream
mathematics is impredicative deserves to be considered further. Indeed,
are there {\it any} compelling examples of clearly mainstream results
that are fundamentally impredicative? And if there are, what are the
philosophical consequences?

First of all,
I do not accept the idea that the existence of even a single impredicative
mainstream theorem would decisively discredit predicativism. If there were
only a handful of examples of impredicative mainstream results and none of
them seemed centrally important, then anyone who found predicativism to be
philosophically persuasive might be willing to simply give up those results.
(Actually, as I will explain in Section 6, this would be an overreaction:
impredicative mathematics does not have to be given up, but only to be
reinterpreted in formalist terms.) But I have already made the point in
Section 3 that the vast bulk of mainstream mathematics is uncontroversially
agreed to be predicative. So the case for predicativism does not hinge on
whether there are {\it any} good examples of clearly mainstream results that
are essentially impredicative.

Having said this, I will now add that I believe there are in fact no
such examples. This runs against claims routinely made in the foundations
literature that certain mainstream results (e.g., Kruskal's theorem)
are known to be impredicative. The basis for assertions like these
is an analysis of the proof-theoretic strength of predicative
systems that was carried out by Feferman and Sch\"utte, based on an idea
of Kreisel. I have thoroughly criticized this analysis, including a series of
subsequent papers of Feferman, in \cite{W3}; here I will outline my critique
of the original analysis involving autonomous systems.

Without going into great detail, the general idea of the autonomous systems
is the following. We construct a recursive well-ordering $\prec$ of $\omega$
with least element $0$ and a corresponding family of formal systems $S_a$
($a \in \omega$) such that we can predicatively infer, for any $a$, that
if all $S_b$ with $b \prec a$ are valid then so is $S_a$. In particular,
the base system $S_0$ is predicatively acceptable; the successor system to
$S_a$ is something like $S_a + {\rm Con}(S_a)$ (actually, it is a bit stronger
than this). Also let ``$a$ is an ordinal notation'' mean that
$\{b \in \omega: b \prec a\}$ is well-ordered.

The question is which systems are predicatively acceptable.
Kreisel's answer was that if $S_0$ proves that $a_1$ is an ordinal notation,
then $S_{a_1}$ can be accepted, and if $S_{a_1}$ proves that $a_2$ is an
ordinal notation, then $S_{a_2}$ can be accepted, and so on.
The Feferman/Sch\"utte analysis identifies a countable ordinal $\Gamma_0$
with the property that it is the smallest ordinal such that there is no
finite sequence of ordinal notations $a_1, \ldots, a_n$ with $a_1 = 0$,
$a_n$ a notation for $\Gamma_0$, and such that $S_{a_i}$ proves that
$a_{i+1}$ is an ordinal notation.  Thus, $\Gamma_0$ is the smallest
predicatively non-provable ordinal. One can show that Kruskal's theorem
implies that a notation for $\Gamma_0$ is well-ordered, so we conclude
that it cannot be predicatively proven.

There are two problems with this analysis. First,
the plausibility of Kreisel's proposal hinges on our conflating two
versions of the concept ``well-ordered'' --- supports transfinite
induction for arbitrary sets versus supports transfinite induction
for arbitrary properties --- which are not predicatively equivalent.

When we prove $a_i$ is an ``ordinal notation'' in the Feferman/Sch\"utte
analysis, we are {\it only} showing that transfinite induction up to $a_i$
holds for statements of the form ``$b$ is in $X$''.  That is, we know that
if the assertion ``for all $b$, everything less than $b$ is in $X$ implies
$b$ is in $X$'' holds for some set $X \subseteq \omega$, then everything less
than $a_i$ is in $X$. {\it But to infer soundness of $S_{a_i}$ from this fact,
we would have to be able to form the set
$$X = \{b: S_b\hbox{ is sound}\},$$
which predicatively we cannot do.} The problem is that the systems $S_a$
are formulated in the language of second order arithmetic and hence
involve assertions about all sets of natural numbers. So in order to
diagnose whether a given $S_a$ is sound we need to quantify over $\Pc(\Nb)$.
This means that the set $X$ displayed above is a set of natural
numbers which is defined by means of a condition that quantifies over
all sets of natural numbers. This is circular in exactly the same way
as the example I described near the end of Section 1. Indeed, it is a
fundamental feature of predicativism that such constructions are considered
illegitimate. The set $X$ is not a predicatively legitimate set, so we
cannot use the fact that transfinite induction up to $a_i$ holds for
sets to predicatively infer that $S_{a_i}$ is sound, and hence the
Feferman/Sch\"utte analysis collapses. This point is made in greater
detail in \cite{W3}. I show there that essentially the same problem
pervades all of the analyses of predicativity which link it to $\Gamma_0$.

I made this criticism public in 2005 and I have yet to see any
substantive counterargument, or even any indication that it has been
clearly understood.

There is a second basic problem. Let $I(a)$ be the assertion that
$a$ is an ordinal notation. In order for $\Gamma_0$ to be the exact
bound for predicativity on the Feferman/Sch\"utte analysis, it must be the
case that for any formula $A$ and any numbers $a$ and $n$ the predicativist
has some way to make the deduction
$${\rm from}\quad I(a)\quad {\rm and}\quad
\hbox{$S_a \vdash A(n)$}\quad
{\rm infer}\quad A(n).\eqno(*)$$
($S_a \vdash A(n)$ means: $A(n)$ is a theorem of $S_a$.)
We have just seen that this inference is predicatively illegitimate on
its face. But suppose the predicativist had some way to draw this inference.
In order for the Feferman/Sch\"utte analysis to succeed, it would also have
to be the case that he {\it cannot} accept the single assertion
$$(\forall a)(\forall n)(I(a) \wedge S_a \vdash A(n)\quad
\to\quad A(n))\eqno(**)$$
for arbitrary $A$, as this would allow him to go beyond $\Gamma_0$.
Not only is ($*$) predicatively invalid, it is hard to imagine
what could lead anyone, predicativist or not, to accept every instance
of ($*$) but not accept ($**$). In \cite{W3} I discuss three separate
remarks of Kreisel on this question; the first is a brief comment that
has no substance, the second is clearly fallacious, and the third is a
highly implausible speculation. As far as I know, no other writer has
even attempted to address this objection, or perhaps even recognized
its existence.

An extra difficulty about this point that has apparently never been
recognized is that whatever reason is given for the predicative
acceptability of each instance of ($*$) would itself have to be impredicative.
For if a predicativist could see that he would accept each instance of ($*$),
then he would grasp the validity of the entire Feferman/Sch\"utte
construction and thus get beyond $\Gamma_0$. Expositions of the
Feferman/Sch\"utte analysis typically presume that the predicative
acceptability of each instance of ($*$) is obvious, which raises
the question of why this is not obvious to the predicativist himself
(as well as the question of why the predicative acceptability of
($**$) is not equally obvious).

I believe this criticism of the Feferman/Sch\"utte analysis is decisive.
Later attempts to verify the same conclusion about $\Gamma_0$ in other ways
involve so many errors (see \cite{W3}) that its uncritical acceptance by the
foundations community for forty years raises serious sociological questions.
There is no connection between predicativism and $\Gamma_0$, nor does there
appear to be {\it any} coherent foundational stance which would lead one to
accept all ordinals less than $\Gamma_0$ and not $\Gamma_0$ itself.

In the second part of \cite{W3} I show how hierarchies of
Tarskian truth predicates
can be used to access ordinals well beyond $\Gamma_0$, sufficient to prove
Kruskal's theorem. I also make a case that these theories are predicatively
legitimate. The basic idea here is to create a hierarchy of formal systems
$(S_a)$ in which the successor theory to $S_a$ has a truth predicate $T_a$
with which one can reason about the truth of statements in the language of
$S_a$. This is similar to the autonomous progressions of Kreisel. But we
then jump one level up and consider a formal system $S_0^1$ in which we can
reason about the entire sequence $(S_a)$, use this as the basis of a new
sequence $(S_a^1)$, and iterate this process. The system $S_0^\omega$
proves that a notation for $\Gamma_0$ is
well-ordered, but there is no reason to stop here. Indeed, we can reason
about the sequence $(S_0^a)$, and this leads to a system which proves that
a notation for the Ackermann ordinal $\phi_{\Omega^2}(0)$ is well-ordered.
This technique can be pushed further to access even larger ordinals.

It is important to
be clear that there are two aspects to \cite{W3}, the first being a
criticism of the various analyses which have been put forward as
showing that $\Gamma_0$ is a limit to predicativism, and the second
being a new analysis involving hierarchies of truth predicates which
goes beyond $\Gamma_0$. If the second part is right, then examples of
impredicative mainstream results vanish almost entirely. The best
example of which I am aware that is not covered by the
results of \cite{W3} is the graph minor theorem; however, it seems
likely that further work along similar lines will succeed in establishing
this result predicatively as well.

Such a large fraction of mainstream mathematics is already uncontroversially
recognized to be predicative that failure to cohere with mathematical
practice is no longer a meaningful criticism of predicativism. If anything,
it is a greater problem for ZFC; see Section 4. Moreover, if the truth theories
approach of \cite{W3} is predicatively valid then meaningful examples
of impredicative mainstream mathematics may disappear altogether.

\section{Set theory and formalism}

ZFC is not a good choice to be the standard foundation for mathematics.
It is unsuitable in two ways. Philosophically, it makes sense only
in terms of a vague belief in some sort of mystical universe of
sets which is supposed to exist aphysically and atemporally (yet,
in order to avoid the classical paradoxes, is somehow ``not there
all at once''). Pragmatically, ZFC fits very badly
with actual mathematical practice insofar as it postulates a vast
realm of set-theoretic pathology which has no relevance to mainstream
mathematics. We might say that it is both theoretically and practically
unsuited to the foundational role in which it is currently cast.

As we have seen, these two defects are linked: reducing to the philosophically
more sensible attitude of considering structures which could conceivably
be physically instantiated has the effect of neatly eliminating the
set-theoretic pathology which characterizes ZFC, while retaining all
of the structures that are essential to mainstream mathematics.

However, this is not to say that the study of ZFC must be abandoned.
It can still be understood as an interesting formal system which
some mathematicians may find quite appealing. The same could be
said for ZFC plus various large cardinal axioms, or ZF plus the
axiom of determinacy, or ZFC plus V = L, or ZFC plus Martin's axiom,
and so on.

How are we to interpret theorems proven in such systems? First of all,
the fact that some statement is a theorem of ZFC is a combinatorial fact
which is just as valid
predicatively as it is classically. We can also regard such a theorem
as a true statement in any model of ZFC, which predicatively exist
provided ZFC is consistent. (And even if we are not certain that ZFC is in
fact consistent, we can still reason under this hypothesis.) What we cannot
do is {\it to regard theorems of ZFC as being true statements in some
canonical universal model of ZFC}, because we have no reason to believe
that such a model exists.

It is worth noting that whatever set theorists privately
believe, their professional behavior is quite consistent with
a formalist interpretation of ZFC. Indeed, a large part of
modern set theory explicitly concerns itself with the study of ZFC
and various extensions of ZFC {\it as formal systems}. For example, the
seminal result of modern set theory is the relative independence of
the continuum hypothesis from ZFC. This is the fact, provable in PA
(or even weaker systems), that if ZFC is consistent then so are
ZFC + CH and ZFC + $\neg\,$CH. One finds in the set theory literature
very little discussion of whether questions like CH are {\it really
true}, but a great deal of work --- work which can be formalized in
PA --- about which axioms imply or fail to imply questions like CH.
Thus, should predicative foundations be universally adopted, the actual
practice of set theory would have to change very little.

A parallel can be drawn between set theory and nonstandard analysis.
In both cases we formally legitimize a vague intuitive idealization
--- actual infinitesimals in the case of nonstandard analysis, and
concretely existing uncountable structures in the case of set theory.
Both are very elegant from a formal standpoint, but not from the
point of view of models: the nonstandard reals do not have the simple,
canonical quality that the standard reals possess, and even if one
believes that ZFC has a canonical model this is only an abstract fact;
we do not have a clear picture of it in the same way that we have a
picture of $\Nb$ (the canonical model of PA) or even $J_2$ (which models
a predicative set-theoretic system \cite{W2}).

Adopting a predicative foundational outlook does not entail abandoning
ZFC altogether. Rather, it means that we must interpret ZFC in formalist
terms, something that, operationally, set theorists already do.

\section{ZFC and number theory}

I just made the point that we have no a priori reason to believe that ZFC
has a canonical model. One consequence of this is that we should be skeptical
of the actual truth of number-theoretic results proven in ZFC that are not
provable in predicative systems. These could only be trusted if we had some
reason to believe that ZFC has models in which $\omega$ is standard.

Now we may believe that ZFC is probably
consistent because (1) no inconsistency has been found yet and (2)
we have built up some sort of intuition for ZFC which tells us
that it is consistent. I personally find these arguments persuasive
but not compelling. They suggest that ZFC probably does have a model.
{\it However, they tell us nothing about whether it has a model with a
standard $\omega$.} This seems to me more likely to be false than true.
Given the recursive compexity of ZFC (as measured by its proof-theoretic
ordinal, and already suggested by the circularity of the power set of $\omega$
mentioned in Section 1) we should not expect that there is such a model
absent some special reason to do so. The presumption should be that
ZFC has no such model and hence that there are probably some false
statements of first order arithmetic that are provable in ZFC.

Antiplatonistic belief in the arithmetical validity of ZFC seems to be
mainly a matter of faith. One could argue that the hierarchy of large cardinal
axioms exhibits a compelling structure which is evidence for the truth
of the arithmetical consequences of these axioms. Maybe so, but this is
at best very indirect evidence and hardly seems very convincing. At present
I think a rational assessment of the evidence would have to conclude that
ZFC very likely proves false number-theoretic assertions.

I hasten to add that this is not an indictment of mainstream number
theory, since mainstream number theory can be formalized in predicative
systems. Rather, we should be suspicious of any number theoretic result
whose proof requires set-theoretically substantial mathematics. (Harvey
Friedman has given examples of such results.)

\section{Conclusion}

In brief, my position is as follows:

Reifying collections as ``abstract objects'' is an elementary philosophical
mistake and is directly responsible for the paradoxes of naive set theory.

The ``iterative conception of sets'' does not succeed in legitimizing
abstract set theory. It crucially involves an idea of {\it set formation}
as if sets were physical objects that could be manipulated, which they are
simultaneously denied to be (yet no characterization of the sort of entities
they are is given). It also takes the ability to form power sets as basic,
yet there is no meaningful sense in which one can imagine {\it forming}
power sets of infinite sets.

We cannot evade the problem of justifying set theory by settling for
the fact that ``mathematicians (meaning set theorists) like ZFC'' or the
probability that ZFC is consistent. At a minimum, if ZFC is to be taken
as the standard foundation of mathematics then we should at least demand
some conviction that it is $\Sigma_1$-valid, which anti-platonism is unable
to provide. Thus ZFC has no clear philosophical basis. We should not ignore
the real possibility that some number-theoretic assertions provable in ZFC
might be false.

Sets are universally understood in quasi-physical terms and their
properties are justified in terms of imagined quasi-physical manipulations.
(E.g., the axiom of choice seems true because we could run through any
family of nonempty sets and choose one element from each, etc.) This is only
legitimate for {\it in principle physically possible} structures, which is
precisely the world of predicativism. The idea that uncountable structures
are in any meaningful sense physically possible does not withstand scrutiny;
indeed, by the L\"owenheim-Skolem theorem we know that we have no way
(that does not presuppose set theory) of describing a possible universe
which contains uncountable structures.

All of the mathematics that is currently applied in science is
predicative, and the idea that this will some day change is not credible.
Virtually all mainstream mathematics is predicative \cite{Sim}, and
it is probably the case that absolutely all mainstream mathematics is
predicative \cite{W3}. The predicative system CM of \cite{W5} is
at least comparable to ZFC in terms of elegance and ease of use.

I am calling for the abandonment of ZFC {\it as a foundational standard}.
I believe its erroneous identification as the correct framework for
mathematics as a whole has led it to receive a disproportionate amount
of attention. However, I am not calling for the study of ZFC to cease
altogether. I believe it and various related systems (including its
augmentation by large cardinals, and variants not including the axiom
of choice) are interesting and worthy of study, but they are also
peripheral to the concerns of mainstream mathematics.



\begin{thebibliography}{aaaaaaaa}
\bibitem {Par}
C.\ Parsons, What is the iterative conception of set?, in \cite{BP},
503-529.
\end{thebibliography}
\end{document}